\documentclass[12pt]{article}
\usepackage{theorem}
\newtheorem{lemma}{Lemma}

\newtheorem{theorem}[lemma]{Theorem}
\newtheorem{corollary}[lemma]{Corollary}
{\theorembodyfont{\upshape}}
{\theorembodyfont{\upshape}}
{\theorembodyfont{\upshape}\newtheorem{example}[lemma]{Example}}
{\theorembodyfont{\upshape}}
\newcommand{\R}{{\bf R}}

\newcommand{\rme}{{\rm e}}
\newcommand{\rmd}{{\rm d}}
\newcommand{\sig}{\sigma}
\newcommand{\alp}{\alpha}
\newcommand{\bet}{\beta}
\newcommand{\gam}{\gamma}
\newcommand{\lam}{\lambda}
\newcommand{\del}{\delta}
\newcommand{\eps}{\varepsilon}
\newcommand{\ome}{\omega}

\newcommand{\lap}{{\Delta}}
\newcommand{\gra}{\nabla}

\newcommand{\tr}{{\rm tr}}
\newcommand{\Dom}{{\rm Dom}}

\newcommand{\norm}{\Vert}

\newcommand{\supp}{{\rm supp}}
\newcommand{\Proof}{\underbar{Proof}{\hskip 0.1in}}
\newcommand{\Note}{{\bf Note}{\hskip 0.1in}}

\newcommand{\dist}{{\rm dist}}

\newcommand{\la}{{\langle}}
\newcommand{\ra}{{\rangle}}
\newcommand{\pr}{\prime}
\newcommand{\nolabel}{\nonumber}
\title{SHARP BOUNDARY ESTIMATES \\
FOR ELLIPTIC OPERATORS}
\author{E.B. Davies}
\date{August 1998}
\begin{document}
\maketitle
%%%%%%%%%%%%%%%%%%%%%
%%%%%%%%%%%%%%%%%%%%%
%\end{document}
%%%%%%%%%%%%%%%%%%%%%
\begin{abstract}
We prove sharp $L^{2}$ boundary decay estimates for the 
eigenfunctions of certain second order elliptic operators 
acting in a bounded region, and of their first space 
derivatives, using only the Hardy 
inequality. These imply $L^{2}$ boundary decay properties of 
the heat kernel and spectral density. We 
deduce bounds on the rate of convergence of the 
eigenvalues when the region is slightly reduced in size. It 
is remarkable that several of the bounds do not involve the space 
dimension.
\par
keywords: boundary decay, Laplacian, Hardy 
inequality, eigenfunctions, heat kernel, spectral density, 
spectral convergence.
\end{abstract}
%%%%%%%%%%%%%%%%%%%%%%%%
\section{Introduction}
\par
Let $H$ be a non-negative second order elliptic operator acting in 
$L^{2}(U,\rmd^{N} x)$ subject to Dirichlet boundary conditions, where 
$U$ is a bounded region in $\R^{N}$ or even in a Riemannian 
manifold. Let $d$ be 
a continuous function on $U$ satisfying $|\gra d|\leq 1$, 
for example the distance from the boundary of $U$, 
which may be very irregular. We say 
that $H$ satisfies Hardy's inequality with respect to $d$ if 
\[
\int_{U}\frac{|f|^{2}}{d^{2}} \leq c^{2}\left( Q(f) +a\norm f 
\norm^{2 }\right)
\]
for all $f\in C_{c}^{\infty}(U)$, where $Q$ is the quadratic 
form of $H$. The precise value of the 
constant $c$ will be of great importance below, but the 
size of $a$ is not crucial.

We are concerned with boundary decay of the eigenfunctions of 
$H$, and more generally of any functions in the domain of 
$H$. Such bounds were first obtained in \cite{EHK, D1} by an 
argument related to that which we use below. The bounds were 
used in \cite{D1} to obtain explicit upper bounds on
the rate at which 
\[
|\lam_{n}(U)-\lam_{n}(U_{\eps})|
\]
vanishes as $\eps\to 0$, where $\lam_{n}(S)$ denotes the $n$-th 
Dirichlet eigenvalue of any region $S$ and 
\[
U_{\eps}:=\{ x\in U:d(x)>\eps \}.
\]
In two recent papers Pang \cite{P1,P2} used a  
different method to obtain a sharp rate  of convergence of 
the first eigenvalue as $\eps \to 0$ for a Dirichlet 
Laplacian in any simply connected subregion of $\R^{2}$, 
and for certain bounded regions in $\R^{N}$. 
In this paper we return to the method of \cite{D1} and show 
how to obtain sharp estimates of the rate of convergence for 
all eigenvalues; see Section 4. 

The key is to obtain better boundary decay estimates of the 
eigenfunctions, in terms of
\[
\int_{\{ x:d(x) <\eps\} } |f|^{2}
\]
for all $f\in\Dom (H)$ and all $\eps >0$, instead of 
estimating
\[
\int_{U}\frac{|f|^{2}}{d^{\gam}}
\]
for all possible $\gam >0$. It is well known to harmonic 
analysts that the former type of estimate is generally 
sharper than the latter, and we find that it yields the 
optimal power in the subsequent 
proof of the convergence of the eigenvalues.
\par
The main theorems of the paper in Section 3 apply to 
weighted Laplace-Beltrami operators acting in regions with 
irregular and possibly fractal boundaries, but in Section 5 
we show that the methods 
can be applied to second order uniformly elliptic  
operators with measurable highest order coefficients. The 
estimates are proved for functions in the domains of the 
operators, and apply in particular to eigenfunctions. In 
most theorems we prove that we have the optimal power of 
$\eps$ in the estimates.

The methods which we use do not require $U$ to be a region 
in a Riemannian manifold. If $U$ is a piecewise manifold 
obtained by glueing together manifolds of the same dimension 
along certain common edges, the same ideas can be applied 
provided the operator $H$ is defined by means of the 
appropriate quadratic form.
\par
In Section 6 we use the results to obtain some new $L^{2}$ 
boundary decay estimates for the heat kernel of the 
operator, and remark that the same methods can be used for 
the spectral density.
\par
The sharp constant in Hardy's inequality is the only 
important input to the argument, and we refer to \cite{D4} for a 
recent review of this topic. Here we mention only a 
few outstanding results for $H:=-\lap_{DIR}$ acting in 
a bounded region $U$ in Euclidean space. If $U$ is a simply connected 
proper subregion of 
$\R^{2}$ then Hardy's inequality holds with 
$c=4$ and $a=0$ by \cite{A}, \cite[Th. 1.5.10]{HKST}. If $U$ is a 
convex region in $\R^{N}$ then it is a matter of folklore 
that Hardy's inequality holds with $c=2$ and 
$a=0$; some proofs are described in \cite{D4}. Finally, if 
$U$ has smooth boundary 
then Hardy's inequality holds with $c=2$ for some $a<\infty$, \cite{BM}.
%%%%%%%%%%%%%%%%%%%%%%%%%%%%%%%%%
\section{Definitions}

We follow the notation of \cite{D1} but with somewhat more 
restrictive conditions on the various coefficient 
functions. Let $\sig$ be a measurable function on 
the incomplete Riemannian manifold $U$ which is positive 
almost everywhere and locally $L^{2}$ wth respect to the 
Riemannian volume element. Define the weighted space $L^{2}(U)$ to be 
the set of (equivalence classes modulo null sets of) functions such that 
\[
\norm f\norm_{2}^{2}:=\int_{U}|f|^{2}\sig^{2} <\infty.
\]
This and subsequent integrals are evaluated using the 
Riemannian volume element.
The introduction of the weight $\sig$  
allows extra applications of our theorems at no cost.
Let $V$ be a non-negative locally $L^{1}$ function on $U$ 
and let $H$ be the operator on $L^{2}(U)$ defined formally by
\[
Hf:=-\sig^{-2}\gra \cdot \left( \sig^{2} \gra f \right) +Vf
\]
subject to Dirichlet boundary conditions. Rigorously we 
start from the non-negative quadratic form
\[
Q(f):=\int_{U}\left( |\gra f|^{2} +V|f|^{2}\right)\sig^{2} 
\]
which is well-defined on the domain $C_{c}^{\infty}(U)$ by 
the hypothesis on $\sig$. We assume that $Q$ is closable 
and define $H$ to be the self-adjoint operator 
on $L^{2}(U)$ associated with the closure 
of the form as described in \cite[Ch. 4]{OPS} and 
\cite[Section 1.2]{HKST}. For discussions of conditions on $\sig$ 
which imply that $Q$ is closable see \cite{RW} and 
\cite[Section 1.2]{HKST}.

If $U$ is a region in $\R^{N}$, $\sig =1$, $V=0$ and 
we choose the Euclidean metric 
then $H=-\lap$ subject to Dirichlet boundary conditions. 
Our results are new in this case when $U$ is bounded and its boundary 
$\partial U$ is irregular, possibly fractal, improving on 
the recent theorems in \cite{D1,P1,P2}. 

Our main assumption is formulated in terms of
a positive continuous function $d$ on $U$ such 
that $|\gra d|\leq 1$, in the weak sense. More precisely, we 
assume that
\[
|d(x)-d(y)|\leq |x-y|
\]
for all $x,y\in U$. This is equivalent to the statement 
that $d$ has distributional 
derivative $\gra d\in L^{\infty}$ which satisfies 
$|\gra d (x)|\leq 1$ almost everywhere in $U$.
One might take $d(x)$ to be the distance 
of $x\in U$ from any closed subset $S$ of $\partial U$ or 
from a closed subset of $M\backslash U$ if $U$ is embedded 
in some larger Riemannian manifold $M$. We assume 
throughout the paper that for 
some constant $c\geq 2$ and some non-negative constant 
$a$ the Hardy inequality (HI)
\[
\int_{U}\frac{|f|^{2}}{d^{2}} \sig^{2}\leq c^{2}\left( Q(f) +a\norm f 
\norm^{2 }\right)
\]
is valid for all $f\in C_{c}^{\infty}(U)$, and hence for all $f$ in 
the domain of the closure of $Q$. Our goal is 
to obtain a similar but stronger bound for all $f\in \Dom 
(H)$ and hence for all eigenfunctions of $H$. Note  
that we do not assume that $U$ is bounded or $H$ has discrete 
spectrum.

A serious difficulty is the fact that we cannot identify 
the domain of $H$ with any of the Sobolev or other 
spaces. If the coefficients of $H$ or the boundary are 
irregular the 
operator domain of $H$ changes if we vary $\sig$ or $V$ within 
the permitted classes, even though the quadratic form domain 
may be unchanged. The bounds which we obtain in Theorems 4 
and 7 bear some relationship with Morrey space 
estimates, already known to be of great importance in the 
theory of elliptic operators, \cite{Au,Gi,Mo}.
%%%%%%%%%%%%%%%%%%%%%%%%%%%%%%%%%%
\section{The main theorems}

Our estimates involve a positive parameter $\eps$, and 
various other constants which depend only on $c\geq 2$, in 
a way which we make explicit. Given $\eps >0$ we put
\[
\ome(x):=\left( \max\{ d(x),\eps\}\right)^{-1/c}
\]
for all $x\in U$.
%%%%%%%%%%%%%%%%%%%%%%%%%%%%%%%%%%%
\begin{lemma}
If $f\in \Dom (H)$ and $s\geq 0$ then
\[
\left|\la Hf,\ome^{2}f\ra +s\norm \ome f\norm_{2}^{2}\right| \leq
c^{2/c}\norm (H+s)f\norm_{2}\norm (H+a)^{1/c}f\norm_{2}.
\]
\end{lemma} 

\Proof Using HI and \cite[Lemma 4.20]{OPS} we have
\[
\ome^{4}\leq (d^{-2})^{2/c}\leq \{c^{2}(H+a)\}^{2/c}
\]
so
\[
0\leq (H+a)^{-1/c}\ome^{4}(H+a)^{-1/c}\leq c^{4/c}I
\]
and
\[
\norm \ome^{2}(H+a)^{-1/c}\norm\leq c^{2/c}.
\]
Hence
\begin{eqnarray*}
	  & &  \left|\la Hf,\ome^{2}f\ra +s\norm \ome f\norm_{2}^{2}\right| \\
	 && = \left| \la (H+s)f,\ome^{2}f\ra \right|   \\
	 && = \left| \la (H+s)f,\ome^{2}(H+a)^{-1/c}.(H+a)^{1/c}f\ra \right|   \\
 & &\leq  c^{2/c} \norm (H+s)f\norm_{2}\norm (H+a)^{1/c}f\norm_{2}. 
\end{eqnarray*}
%%%%%%%%%%%%%%%%%%%%%%%%%%%%%%%%%%%
\begin{lemma}
If $f\in \Dom(Q)$ and $\mu\in W^{1,\infty}(U)$ then $\mu f \in \Dom (Q)$ and
\[
Q(\mu f)\leq 2\norm\mu\norm_{\infty}^{2}Q(f)
+2\norm\gra\mu\norm_{\infty}^{2}\norm f\norm_{2}^{2}.
\]
\end{lemma}

\Proof If $f\in C_{c}^{\infty} (U)$ then $\mu f\in 
W^{1,\infty}_{c}(U)\subseteq \Dom(Q)$ and
\begin{eqnarray*}
	Q(\mu f) & = & \int_{U}\left( |\mu \gra f +f\gra 
	\mu|^{2}+V\mu^{2}|f|^{2}\right)\sig^{2}  \\
	 & \leq & \int_{U}\left(2\mu^{2}|\gra f|^{2} +2|f|^{2}|\gra 
	 \mu|^{2}
	 +V\mu^{2}|f|^{2}\right)\sig^{2}\\
	 & \leq & 2\norm 
	 \mu\norm_{\infty}^{2}Q(f)+2\norm\gra\mu\norm_{\infty}^{2}\norm 
	 f\norm_{2}^{2}.
\end{eqnarray*}
If $f\in\Dom(Q)$ then the fact that $\mu f\in\Dom(Q)$ and 
the validity of the same estimate both follow from the lower 
semi-continuity of $Q$.
%%%%%%%%%%%%%%%%%%%%%%%%%%%%%%
\begin{lemma}
If $f\in\Dom(H)$ then
\[
\int_{U}\frac{\ome^{2}|f|^{2}\sig^{2}}{c^{2}d^{2}} \leq 
c^{2/c}\norm (H+a)f\norm_{2}\norm 
(H+a)^{1/c}f\norm_{2}+\int_{U}|\gra \ome |^{2}|f|^{2}\sig^{2}.
\]
\end{lemma}

\Proof By HI and Lemma 2 we see that $\ome f\in\Dom(Q)$ and
\[
\int_{U}\frac{\ome^{2}|f|^{2}\sig^{2}}{c^{2}d^{2}} \leq 
Q(\ome f)+a\norm \ome f\norm_{2}^{2}.
\]
Secondly
\begin{eqnarray}
	 &  & Q(\ome f)-\frac{1}{2}\la Hf,\ome^{2}f\ra 
	 - \frac{1}{2}\la \ome^{2}f,Hf\ra \nolabel\\
	 &  & =\int_{U}\left\{ |\gra (\ome f)|^{2}- 
	 \frac{1}{2}\gra f\cdot\overline{\gra (\ome^{2}f)} 
	 -\frac{1}{2}\gra (\ome^{2}f)\cdot \overline{\gra f} 
	 \right\}\sig^{2} \nolabel\\
	 &  & =\int_{U}|\gra \ome |^{2}|f|^{2}\sig^{2}.  
	 \label{alpha}
\end{eqnarray}
Hence
\begin{eqnarray*}
	 &  & Q(\ome f)+a\norm \ome f\norm_{2}^{2}  \\
	 &  & =\frac{1}{2}\left( \la Hf,\ome^{2}f\ra +a\norm \ome 
	 f\norm_{2}^{2}\right)
	+ \frac{1}{2}\left( \la \ome^{2}f,Hf\ra+a\norm \ome 
	 f\norm_{2}^{2}\right) +\int_{U}|\gra \ome |^{2}|f|^{2}\sig^{2}.
\end{eqnarray*}
The proof is completed by combining the above two formulae 
with the bound of Lemma 1.

For some comments on the optimality of the estimates in the 
following theorem see Example \ref{A1}, Example \ref{B1} and 
the note after Corollary \ref{C1}.
%%%%%%%%%%%%%%%%%%%%%%%%%%%%%%%%%%%%%
\begin{theorem}\label{four}
If $f\in\Dom(H)$ then assuming HI we have
\[
\int_{\{x:d(x) <\eps  \}}\frac{|f|^{2}}{d^{2}}\sig^{2}\leq 
c_{0}\eps^{2/c}\norm (H+a)f\norm_{2} 
\norm(H+a)^{1/c}f\norm_{2}
\]
for all $\eps >0$, where $c_{0}:=c^{2+2/c}$. Hence
\[
\int_{\{x:d(x) <\eps  \}}|f|^{2}\sig^{2}\leq 
c_{0}\eps^{2+2/c}\norm (H+a)f\norm_{2} 
\norm(H+a)^{1/c}f\norm_{2}
\]
for all $\eps >0$.
\end{theorem}

\Proof We rewrite Lemma 3 in the form
\[
\int_{U}Y|f|^{2}\sig^{2}\leq c^{2/c}\norm (H+a)f\norm_{2} 
\norm(H+a)^{1/c}f\norm_{2}
\]
where
\[
Y:=\frac{\ome^{2}}{c^{2}d^{2}} -|\gra \ome |^{2}.
\]
If $d(x)\geq \eps$ then 
\[
|\gra\ome |^{2} \leq c^{-2}d^{-2-2/c}=\frac{\ome^{2}}{c^{2}d^{2}}
\]
so $Y(x)\geq 0$. On the other hand if $d(x)<\eps$ then
\[
Y=\frac{\ome^{2}}{c^{2}d^{2}}\geq 
\frac{1}{c^{2}\eps^{2/c}d^{2}}.
\]
Thus 
\begin{eqnarray*}
	\int_{\{x:d(x) <\eps  \}}\frac{|f|^{2}}{d^{2}}\sig^{2} &\leq & 
	c^{2}\eps^{2/c}\int_{U}Y|f|^{2}\sig^{2}  \\
	 & \leq & c^{2+2/c}\eps^{2/c}\norm (H+a)f\norm_{2} 
\norm(H+a)^{1/c}f\norm_{2}.
\end{eqnarray*}
The second statement of the theorem is an immediate 
consequence of the first.
%%%%%%%%%%%%%%%%%%%%%%%%%%%%%%%%
\begin{example} \label{A1} If $c=2$ and $\sig=1$ then the theorem states that
\[
\int_{\{x:d(x) <\eps  \}}|f|^{2}\leq 
8\eps^{3}\norm (H+a)f\norm_{2} 
\norm(H+a)^{1/2}f\norm_{2}
\]
for all $f\in\Dom(H)$, which is the $L^{2}$ analogue 
of $f(x)=O(d(x))$ as $d(x)\to 0$. In particular suppose that  
$U\subseteq \R^{N}$ is bounded with a smooth boundary 
$\partial U$, and let $f$ be a generic function in 
$C_{c}^{\infty}(\overline{U})$ which vanishes on $\partial 
U$. Then $f\in\Dom (H)$ and the power $3$ of $\eps$ above is optimal.
\end{example}
%%%%%%%%%%%%%%%%%%%%%%%%%%%%%%%%
\begin{example} \label{B1} 
Let $U:=(0,\infty)$, $d(x):=x$ and $\sig(x):=x^{\alp/2}$ where 
$0\leq \alp <1$. Then 
the operator $H$ is given formally by
\[
Hf(x):=-x^{-\alp}\frac{\rmd}{\rmd x} \left\{ x^{\alp}\frac{\rmd 
f}{\rmd x} \right\}
\]
subject to Dirichlet boundary conditions at $x=0$. The 
quadratic form $Q$ has domain $W^{1,2}_{0}((0,\infty), 
x^{\alp} \rmd x)$. A standard result, \cite[p. 104]{STDO}, states that 
the strong Hardy inequality holds with $c=2/(1-\alp)$.

Now let $f$ be a smooth function on $(0,\infty)$ which 
vanishes  for $x>2$ and equals $x^{1-\alp}$ for $0<\alp <1$. 
It is easy to prove that $f\in\Dom (H)\subseteq 
W^{1,2}_{0}((0,\infty),x^{\alp}\rmd x)$. If $0<\eps <1$ then 
one also has
\[
\int_{\{ x:d(x)<\eps\}}\frac{|f|^{2}}{d^{2}}\sig^{2}
=\frac{\eps^{3-\alp}}{3-\alp}=k\eps^{2+2/c}.
\]
Therefore the power $2+2/c$ in Theorem \ref{four} is optimal.
\end{example}
%%%%%%%%%%%%%%%%%%%%%%%%%%%%%%
\begin{corollary}
If $g(s)$ is a monotonically decreasing $C^{1}$ function 
on $(0,\del]$ which vanishes for $s=\del$ then
\begin{equation}
\int_{U}g(d)|f|^{2}\sig^{2}\leq 
c_{0}\int_{0}^{\del}|g^{\pr}(s)|s^{2+2/c}\rmd s . \norm 
(H+a)f\norm_{2}\norm (H+a)^{1/c}f\norm_{2}    \label{gest}
\end{equation}
for all $f\in\Dom(H)$, provided the integral on the RHS 
is finite.
\end{corollary}

\Proof We have
\[
g(d)=\int_{d}^{\del}|g^{\pr}(s)|\rmd s
\]
for all $d\in (0,\del]$. The corollary follows by applying Fubini's 
theorem to
\[
\int_{s=0}^{\del} \left\{ 
\int_{U}\chi_{d<s}|g^{\pr}(s)|\, |f|^{2}\sig^{2}
\right\} \rmd s 
\]
where $\chi$ stands for the characteristic function of a set.

\Note Let $H:=-\lap_{DIR}$ acting in $L^{2}(U,\rmd x)$ 
with $\sig:=1$, where $U$ is a bounded region in $\R^{N}$, 
and let $d$ be the distance to the boundary $\partial U$.
If $g(s)=o(s^{-2-2/c})$ as $s\to 0$ then (\ref{gest}) 
is equivalent to
\begin{equation}
\int_{U}g(d)|f|^{2}\leq 
c_{0}(2+2/c)\int_{0}^{\del}g(s)s^{1+2/c}\rmd s . \norm 
(H+a)f\norm_{2}\norm (H+a)^{1/c}f\norm_{2}    \label{gest2}
\end{equation}
which may be compared with the pointwise bound
\begin{equation}
|f(x)|=O(d(x)^{1/2+1/c}) \label{pointwise}
\end{equation}
as $d(x)\to 0$. If $\partial U$ is $C^{2}$ then $c=2$ and this pointwise 
bound is sharp for the first eigenfunction of $H$. However, 
no such pointwise bound exists for arbitrary functions in 
the domain of $H$. Moreover, if $\partial U$ is fractal it 
is not clear that (\ref{pointwise}) is the correct pointwise 
analogue of (\ref{gest2}), nor indeed that there is any 
pointwise analogue.

Our next task is to obtain comparable estimates for $|\gra 
f|$. This necessitates introducing the continuous function on 
$\tau:U\to [0,\infty)$ defined by
\[
\tau(x):=\left\{ \begin{array}{ll}
\eps^{-1/c}&\mbox{ if $0<d(x)\leq \eps$}\\
c^{-1}\eps^{-1-1/c}((1+c)\eps-d(x))&\mbox{ if $\eps < 
d(x)\leq (1+c)\eps$}\\
0&\mbox{otherwise.}
\end{array}\right.
\]
It is immediate that 
\begin{eqnarray*}
	 0\leq \tau & \leq & \ome \leq \eps^{-1/c}  \\
	|\gra\tau | & \leq & c^{-1}\eps^{-1-1/c}  \\
	\supp(\tau) & \subseteq & \{x: 0\leq d\leq (1+c)\eps\}.
\end{eqnarray*}
%%%%%%%%%%%%%%%%%%%%%%%%%%
\begin{theorem} If $f\in\Dom (H)$ then assuming HI we have
\[
\int_{\{x:d(x) <\eps\}}|\gra f|^{2}\sig^{2}
\leq c_{1}\eps^{2/c}\norm 
(H+a)f\norm_{2}\norm(H+a)^{1/c}f\norm_{2}.
\]
for all $\eps >0$, where
\[
c_{1}:=c^{2/c}+c^{2/c}(1+c)^{2+2/c}.
\]
\end{theorem}

\Proof We have
\begin{eqnarray*}
     &	&\eps^{-2/c}\int_{\{ x:d(x) <\eps\}}|\gra f|^{2}\sig^{2} \\
	  &  &=  \int_{\{ x:d(x) <\eps\}}|\gra (\tau f)|^{2}\sig^{2}  \\
	  &   &\leq Q(\tau f) 
\end{eqnarray*}
where $\tau f\in\Dom(Q)$ by Lemma 2. By the same argument 
as in (\ref{alpha}) of Lemma 3 this equals
\begin{eqnarray*}
	 &  & \frac{1}{2}\la Hf,\tau^{2}f\ra +\frac{1}{2}\la \tau^{2}f,Hf\ra 
	 +\int_{U} |\gra \tau |^{2}|f|^{2} \sig^{2}  \\
	 &  & \leq \norm Hf\norm_{2}\norm 
	 \tau^{2}(H+a)^{-1/c}\norm\, \norm (H+a)^{1/c}f\norm_{2} 
	 + \int_{U} |\gra \tau  |^{2}|f|^{2} \sig^{2} \\
	 &  & \leq c^{2/c}\norm Hf\norm_{2}\norm  (H+a)^{1/c}f\norm_{2} 
	 +c^{-2}\eps^{-2-2/c}\int_{\{x:d(x)<(1+c)\eps\}} |f|^{2}\sig^{2}\\
	 &  & 
	 \leq c_{1} \norm (H+a)f\norm_{2}\norm  (H+a)^{1/c}f\norm_{2} 
\end{eqnarray*}
using Theorem 4.

\Note By extending the calculation of Example \ref{C1} one 
sees that the power of $\eps$ in the above 
theorem is optimal. The choice of $\tau$ in the proof is 
certainly not optimal, so neither is the value of $c_{1}$ 
obtained.

%%%%%%%%%%%%%%%%%%%%%%%%%%%%
\begin{corollary}\label{C1}
If $Hf=\lam f$ and $\norm f \norm_{2}=1$ then
\begin{equation}
\int_{\{x: d(x)<\eps\}}|f|^{2}\sig^{2} \leq 
c_{0}\eps^{2+2/c}(\lam +a)^{1+1/c}\label{eigen}
\end{equation}
and 
\[
\int_{\{x: d(x)<\eps\}}|\gra f|^{2}\sig^{2} \leq 
c_{1}\eps^{2/c}(\lam +a)^{1+1/c}
\]
for all $\eps >0$.
\end{corollary} 

\Proof These follow directly from Theorems 4 and 6.

%%%%%%%%%%%%%%%
\Note If we insert the eigenfunction $f$ directly 
into HI we obtain
\[
\int_{\{x: d(x)<\eps\}}|f|^{2}\sig^{2} \leq 
c^{2}\eps^{2}(\lam +a)
\]
which is exactly what is obtained by interpolating between 
(\ref{eigen}) and the trivial estimate
\[
\int_{\{x: d(x)<\eps\}}| f|^{2}\sig^{2} \leq  1.
\]
This supports the conjecture that the constant $c_{0}$ in 
Theorem \ref{four} is optimal.
%%%%%%%%%%%%%%%%%%%%%%%%%%%%%%%%%
\begin{corollary}
If $H:=-\lap_{DIR}$ in $L^{2}(U,\rmd^{2} x)$ where $U$ is a 
simply connected proper subregion of $\R^{2}$ and 
\[
d(x):=\dist(x,\partial U)
\]
then 
\[
\int_{\{x:d(x) <\eps  \}}|f|^{2}\leq 
32\eps^{5/2}\norm Hf\norm_{2} 
\norm H^{1/4}f\norm_{2}
\]
and
\[
\int_{\{x:d(x) <\eps  \}}|\gra f|^{2}\leq 
114\eps^{1/2}\norm Hf\norm_{2} 
\norm H^{1/4}f\norm_{2}
\]
for all $f\in\Dom(H)$.
\end{corollary}

\Proof We may put $c=4$, $a=0$ and $\sig=1$ in 
Theorems 4 and 6 by \cite{A}, \cite[Th. 1.5.10]{HKST}.

%%%%%%%%%%%%%%%%%%%%%%%%%%%%%%%%%%%%%%%%%
\section{Perturbation of the domain}

In this section we use the results above to consider the 
effect on the spectrum of $H$ of replacing the region $U$ 
by a slightly smaller region $U_{\eps}$ such that
\[
\{ x\in U:d(x) >\eps\} \subseteq U_{\eps}\subseteq U.
\]
If $\lam_{n}(U_{\eps})$ denote the eigenvalues of the 
operator $H_{\eps}$ defined by restricting $H$ to $L^{2}(U_{\eps})$ 
where we again impose Dirichlet boundary conditions, then 
variational arguments imply that 
\[
\lam_{n}(U)\leq \lam_{n}(U_{\eps})
\]
for all $n$ and $\eps >0$, and our goal is to find 
quantitative estimates of the difference. The constants 
$c_{n}$ below all depend only on  $a,c,c_{0},c_{1}$ and $n$. 

Let $\mu:U\to [0,\infty)$ be defined by
\[
\mu(x):=\left\{ \begin{array}{ll}
0&\mbox{ if $0<d(x)\leq \eps$}\\
(d(x)-\eps)/\eps &\mbox{ if $\eps < d(x)\leq 2\eps$}\\
1 &\mbox{otherwise.}
\end{array}\right.
\]
so that $0\leq \mu\leq 1$, $|\gra\mu |\leq \eps^{-1}$ and 
$\mu$ has support in $U_{\eps}$.

%%%%%%%%%%%%%%%%%%%%%%%%%%%%%%%
\begin{lemma} 
There exists a constant $c_{2}\geq 0$ such that if 
$f\in\Dom(H)$ then
\[
Q(\mu f) \leq Q(f) +\eps^{2/c}c_{2}\norm (H+a)f 
\norm_{2}\norm (H+a)^{1/c}f\norm_{2}.
\]
\end{lemma}

\Proof  Putting $S:=\{ x:\eps < d(x) <2\eps \}$ we have
\begin{eqnarray*}
	Q(\mu f)-Q(f) & \leq & \int_{S}|\gra (\mu f)|^{2}\sig^{2}  \\
	 & \leq & 2\int_{S}\mu^{2}|\gra 
	 f|^{2}\sig^{2}+2\int_{S}|\gra \mu |^{2}|f|^{2}\sig^{2}  \\
	 & \leq  & 2\int_{S}|\gra 
	 f|^{2}\sig^{2}+2\eps^{-2}\int_{S}|f|^{2}\sig^{2}   \\
	 & \leq & \eps^{2/c}c_{2}\norm (H+a)f 
\norm_{2}\norm (H+a)^{1/c}f\norm_{2} 
\end{eqnarray*}
by Theorems 4 and 6.

It is crucial to the application of our next lemma that 
$0<\eps^{1+1/c} <\eps^{2/c}$
provided $0<\eps <1$, so the error is actually smaller 
than that of Lemma 9 as $\eps\to 0$.

%%%%%%%%%%%%%%%%%%%%%%%%%%%%%%%
\begin{lemma}
There exists a constant $c_{3}\geq 0$ such that if 
$f\in\Dom(H)$ then
\[
\norm f \norm_{2} \geq \norm \mu f\norm_{2}\geq \norm f 
\norm_{2}-c_{3}\eps^{1+1/c}
\norm (H+a)f \norm_{2}^{1/2}\norm (H+a)^{1/c}f\norm_{2}^{1/2}.
\]
\end{lemma}

\Proof The first inequality is elementary. We also have
\begin{eqnarray*}
\left| \norm f \norm_{2} -\norm \mu f\norm_{2} \right|^{2}
     & \leq & \norm f-\mu f\norm_{2}^{2}\\
	 & = & \int_{U}(1-\mu )^{2}|f|^{2}\sig^{2}  \\
	 & \leq & \int_{\{ x:d(x)<2\eps\}}|f|^{2}\sig^{2}\\
	 & = & \eps^{2+2/c} c_{3}^{2}
	 \norm (H+a)f \norm_{2}\norm (H+a)^{1/c}f\norm_{2}
\end{eqnarray*}
by Theorem 4. The second inequality of the lemma follows.

%%%%%%%%%%%%%%%%%%%%%%%%%
The case $n=1$ of the following theorem with the sharp 
power $\eps^{1/2}$ corresponding to $c=4$ 
was already proved for proper simply 
connected subregions of $\R^{2}$ in \cite{P2}, by an 
entirely different method which seems not to extend to 
higher eigenvalues.

%%%%%%%%%%%%%%%%%%%%%%%%%
\begin{theorem}
There exist constants $c_{n}$ for all positive integers 
$n$ such that 
\[
\lam_{n}(U)\leq \lam_{n}(U_{\eps})\leq 
\lam_{n}(U)+c_{n}\eps^{2/c}.
\]
\end{theorem}

\Proof this follows \cite[Th. 22]{D1} closely.
%%%%%%%%%%%%%%%%%%%%%%%%%%%%%%%
\section{Elliptic operators}

In this section we extend the earlier results to second 
order uniformly elliptic operators in divergence form with possibly 
measurable second order coefficients, making use only of the 
ellipticity constant of the operator. Throughout the section 
we put $\sig=1$ and integrate with respect to Lebesgue measure.

Let $U$ be a bounded region in $\R^N$ with $C^2$ boundary 
and let 
\[
d(x):=\dist(x,\partial U)
\]
so that
\[
\int_U \frac{|f|^2}{d^2} \leq 4\int_U (|\gra f|^2+a|f|^2)
\]
for some $a\geq 0$ and all $f \in W_0^{1,2}(U)$. Now let 
\[
Hf(x):=-\sum_{i,j} \frac{\partial}{\partial x_i} 
\left(a^{i,j} (x) \frac{\partial f}{\partial x_j} \right)
\]
subject to Dirichlet boundary conditions in $L^2(U)$, 
where 
\[
1\leq a(x) \leq \alp^2 
\]
for all $x\in U$, and we 
interpret $H$ as a self-adjoint operator using the theory 
of quadratic forms as usual. If we put $\tilde d(x):=\alp^{-
1} d(x)$ then 
\[
\sum_{i,j} a^{i,j} (x) \frac{\partial \tilde d}{\partial x_i}
\frac{\partial \tilde d}{\partial x_j}\leq 1
\]
for all $x\in U$ and
\[
\int_U \frac{|f|^2}{\tilde d^2} \leq 
4\alp^2\left(Q(f)+a\norm f\norm_2^2\right)
\]
for all $f \in W_0^{1,2}(U)$, where $Q$ is the quadratic 
form asociated with $H$.

We are now in a position to apply the theory of the paper 
to the pair $H,\tilde d$ with $c:=2\alp$.

%%%%%%%%%%%%%%%%%%%%%%%%%%%
\begin{theorem}
There exists a constant $c_{0}$ such that if $f\in\Dom(H)$ then
\begin{equation}
\int_{\{x:d(x) <\eps  \}}|f|^{2}\leq 
c_{0}\eps^{2+1/\alp}\norm (H+a)f\norm_{2} 
\norm(H+a)^{1/(2\alp)}f\norm_{2}.                 \label{B}
\end{equation}
\end{theorem}

\begin{theorem} There exists a constant $c_{1}$ such that 
if $f\in\Dom (H)$ then
\begin{equation}
\int_{\{x:d(x) <\eps\}}|\gra f|^{2}
\leq c_{1}\eps^{1/\alp}\norm 
(H+a)f\norm_{2}\norm(H+a)^{1/(2\alp)}f\norm_{2}.  \label{C}
\end{equation}
\end{theorem}

We next suppose that $U_{\eps}$ is a region satisfying 
the same conditions as in Theorem 11, and define 
$\lam_{n}(U_{\eps})$ in a similar manner.

\begin{theorem}
There exist constants $c_{n}$ for all positive integers 
$n$ such that 
\[
\lam_{n}(U)\leq \lam_{n}(U_{\eps})\leq 
\lam_{n}(U)+c_{n}\eps^{1/\alp}.
\]
\end{theorem}

In each case we conjecture that the power of $\eps$ is 
optimal. The three theorems can be proved in two ways. We 
may adapt the proofs of this paper, replacing the weighted 
Laplacian by a more general second order elliptic operator. 
Alternatively, we may apply the theorems of the paper, but 
using a Riemannian metric and weight adapted to the 
choice of the second order coefficients, as described 
in \cite{D1}. Namely if 
$g_{i,j}(x)$ is the matrix inverse to $a^{i,j}(x)$ then the 
Riemannian metric
\[
\sum_{i,j}g_{i,j}(x)\rmd x^{i}\rmd x^{j}
\] 
is Lipschitz equivalent to the Euclidean metric in $U$. 
Indeed the Riemannian distance function is 
bounded between $\alp^{-1}d$ and $d$. If 
also
\[
\sig(x):=\det\left( g_{i,j}(x) \right)^{-1/4}
\]
then
\begin{eqnarray*}
	\int_{U}|\gra f|^{2}\sig^{2} & = & 
	\int_{U}\sum_{i,j}a^{i,j}(x)
	\frac{\partial f}{\partial x_{i} }\frac{\partial 
	\overline{f}}{\partial x_{j} }\rmd^{N}x \\
 \int_{U}| f|^{2}\sig^{2} & = & \int_{U}|f|^{2}\rmd^{N}x
\end{eqnarray*}
where the integrals on the left are with respect to the 
Riemannian measure and $\rmd^{N}x$ is the Lebesgue measure. 
Hence
\begin{eqnarray*}
	Hf & := & -\sig^{2}\gra \cdot \left( \sig^{2}\gra f\right) \\
	 & = & -\sum_{i,j} \frac{\partial}{\partial x_i} 
\left(a^{i,j} (x) \frac{\partial f}{\partial x_j} \right).
\end{eqnarray*}

We finally remark that the above theorems can be localised. 
Suppose $H$ is sub-elliptic but the hypothesis 
$1\leq a(x) \leq \alp^2$ 
holds for $x$ such that  $\dist (x,S)<\bet$, where $\bet 
>0$ and $S$ is some closed subset of $\partial U$. We only 
assume that $\partial U$ is $C^{2}$ in the 
$\bet$-neighbourhood of $S$. If we put 
\[
\tilde{d}(x)=\alp^{-1}\min \{ \dist(x,S), \bet \}
\]
then
\[
\sum_{i,j} a_{i,j} (x) \frac{\partial \tilde d}{\partial x_i}
\frac{\partial \tilde d}{\partial x_j}\leq 1
\]
for all $x\in U$, because the gradient of $\tilde{d}$ 
vanishes outside the $\bet$-neighbourhood of $S$. The proof 
of 
\[
\int_U \frac{|f|^2}{\tilde d^2} \leq 
4\alp^2\left(Q(f)+a\norm f\norm_2^2\right)
\]
for all $f \in W_0^{1,2}(U)$ involves the same arguments as 
in \cite{BM}, concentrating on the region $\{x\in U:\dist 
(x,S)<\bet\}$.
 
%%%%%%%%%%%%%%%%%%%%%%%%%%%%%
\section{Heat kernel and related bounds}

If $K(t,x,y)$ is the heat kernel of a uniformly elliptic 
second order operator $H$ written in divergence form and 
acting in $L^{2}(U, \rmd^{N}x)$ subject to Dirichlet boundary 
conditions, then it is known that 
\begin{equation}
0\leq K(t,x,y)\leq c_{2}t^{-N/2}   \label{A}
\end{equation}
for all $x,y\in U$ and all $t>0$. If $H=-\lap_{DIR}$ then we 
may even take $c_{2}=(4\pi)^{-N/2}$; see \cite{HKST} for an 
account of the relevant heat kernel bounds. We are 
interested in $L^{2}$ boundary decay properties of the heat kernel 
and spectral density which bear some relationship with 
the `intrinsically ultracontractive'
pointwise bounds obtained under much stronger assumptions 
and with much less control on the constants in \cite[Chapter 
4]{HKST}. Throughout this section we assume (\ref{B}), 
(\ref{C}) and (\ref{A}); the constants in our bounds depend 
on these constants and on $N$ in a manner which is 
easy to make explicit.

%%%%%%%%%%%%%%%%%%%%%%%%%
\begin{theorem}
Under the assumptions (\ref{B}), (\ref{C}) and (\ref{A}) we 
have
\begin{equation}
\int_{\{ x:d(x) <\eps\}} K(t,x,y)^{2}\rmd^{N}x \leq 
c_{3}\rme^{at}(\eps^{2}/t)^{1+1/(2\alp)}t^{-N/2}          \label{ker1}
\end{equation}
for all $y\in U$ and all $t>0$. If also $U$ is bounded then
\begin{equation}
\int_{\{ x:d(x) <\eps\}} K(t,x,x)\rmd^{N}x \leq 
c_{4}(\eps^{2}/t)^{1+1/(2\alp)}t^{-N/2}                \label{ker2}
\end{equation}
for all $\eps >0$ and all $0<t\leq 1$.
\end{theorem}
\Proof Denoting the left-hand side of  (\ref{ker1}) by $I$ we have
\[
I=\int_{\{ x:d(x) <\eps\}}|\rme^{-Ht}\del_{y}|^{2}\rmd^{N}x
=\int_{\{ x:d(x) <\eps\}}|\rme^{-Ht/2}g|^{2}\rmd^{N}x
\]
where $\del_{y}$ is the delta function at $y\in U$ and 
$g:=\rme^{-Ht/2}\del_{y}$. Using the semigroup property we 
have
\[
\norm g \norm_{2}^{2} =\int_{U}K(t/2,x,y)^{2}\rmd^{N}x
=K(t,y,y)\leq c_{2}t^{-N/2}.
\]
Applying (\ref{B}) and then the spectral theorem we deduce
\begin{eqnarray*}
	I & \leq & c_{0}\eps^{2+1/\alp}\norm 
	(H+a)\rme^{-Ht/2}g\norm_{2} \norm 
	(H+a)^{1/(2\alp)}\rme^{-Ht/2}g\norm_{2} \\
	 & \leq & c_{0}\eps^{2+1/\alp}\rme^{at}\norm 
	(H+a)\rme^{-(H+a)t/2}\norm \,\norm 
	(H+a)^{1/(2\alp)}\rme^{-(H+a)t/2}\norm \, \norm g\norm_{2}^{2} \\
	 & \leq & c_{3}\eps^{2+1/\alp}\rme^{at}t^{-1-1/(2\alp)-N/2}.
\end{eqnarray*}

We adopt an alternative strategy to prove the second 
inequality. Let $\{ \lam_{n}\}_{n=1}^{\infty}$ be the 
eigenvalues of $H$ written in increasing order and repeated 
according to multiplicity, and let $\{ 
\phi_{n}\}_{n=1}^{\infty}$ be the corresponding normalised 
eigenfunctions. It is known that there exist positive 
constants $a_{1}$ and $a_{2}$ such that
\[
a_{1}n^{2/N}\leq \lam_{n}\leq a_{2}n^{2/N}
\]
for all $n$. Also 
\[
K(t,x,x)=\sum_{n=1}^{\infty}\rme^{-\lam_{n}t} 
|\phi_{n}(x)|^{2}
\]
for all $x\in U$ and $t>0$. Denoting the left-hand side of 
(\ref{ker2}) by $J$ we deduce that for $0<t\leq 1$
\begin{eqnarray*}
	J & = & \sum_{n=1}^{\infty}\rme^{-\lam_{n}t}
	\int_{\{ x:d(x) <\eps\}} |\phi_{n}|^{2} \rmd^{N}x\\
	 & \leq & \sum_{n=1}^{\infty}\rme^{-\lam_{n}t}
	 c_{0}\eps^{2+1/\alp}(\lam_{n}+a)^{1+1/(2\alp)}  \\
	 & \leq & c_{0}\eps^{2+1/\alp} \sum_{n=1}^{\infty}
	  \rme^{-a_{1}n^{2/N}t}(a_{2}n^{2/N}+a)^{1+1/(2\alp)}\\
	 & \leq & c_{4}\eps^{2+1/\alp}t^{-1-1/(2\alp)-N/2}
\end{eqnarray*}
where for the last line we compared the sum with the 
corresponding integral.

The estimate (\ref{ker2}) is not asymptotically optimal 
as $\eps,t\to 0$, even for $H:=-\lap_{DIR}$ acting in 
$L^{2}(0,\infty )$. In this case we have
\begin{eqnarray*}
	K(t,x,x) & = &   (4\pi t)^{-1/2}\left( 1-\rme^{-x^{2}/t}\right)\\
	 & \sim & (4\pi )^{-1/2}x^{2}t^{-3/2}
\end{eqnarray*}
if $0<x^{2}<<t$ by the reflection principle \cite[p107]{HKST}. 
Therefore
\[
\int_{0}^{\eps}K(t,x,x)\rmd x \sim (36\pi 
)^{-1/2}\eps^{3}t^{-3/2}
\]
if $0<\eps^{2}<<t$. However Theorem 15 with $\alp =1$ only 
yields
\[
\int_{0}^{\eps}K(t,x,x)\rmd x \leq c_{4}\eps^{3}t^{-2}.
\]

The following is a possible reason for this failure. 
Optimal estimates on eigenfunctions associated 
with highly degenerate eigenvalues can be much worse than 
one expects for typical eigenvalues. Our proof uses a 
bound for every eigenfunction which takes no account of this 
fact, so when summed up it is not 
surprising that the resulting heat kernel bound is not optimal.

We note that the same method may be used to obtain upper 
bounds on 
\[
\int_{\{ x:d(x) <\eps\}} 
\sum_{n=1}^{\infty}\rme^{-\lam_{n}t}|\gra \phi_{n}(x)|^{2}\rmd^{N}x. 
\]

Another approach may be used to obtain upper bounds on 
quantities associated with the spectral density. Let 
$E_{\lam}$ be the spectral projection of $H$ associated with 
the interval $(-\infty,\lam)$, and let $e(\lam,x,y)$ be its 
integral kernel. Then
\[
N(\lam)=\int_{U}e(\lam ,x,x)\rmd x
\]
where $N(\lam)$ is the number of eigenvalues of $H$ less 
than $\lam$. If $\eps >0$ we put
\[
N(\eps,\lam):=\int_{\{x:d(x) <\eps  \}}e(\lam,x,x)\rmd x.
\]

\begin{theorem}
Under the assumption (\ref{B}) we have
\[
N(\eps,\lam)\leq c_{0}\eps^{2+1/\alp}(\lam 
+a)^{1+1/(2\alp)} N(\lam)
\]
for all $\lam\geq 0$ and $\eps >0$.
\end{theorem}

\Proof If $f=E_{\lam} f$ then (\ref{B}) implies
\[
\int_{\{x:d(x) <\eps  \}}|f|^{2}\leq 
c_{0}\eps^{2+1/\alp}(\lam+a)^{1+1/(2\alp)}\norm f\norm_{2}^{2}        
\]
which is equivalent to the operator inequality
\[
\norm Q_{\eps}E_{\lam}\norm ^{2}\leq c_{0}\eps^{2+1/\alp}(\lam+a)^{1+1/(2\alp)}
\]
where $Q_{\eps}$ is the projection given by multiplying by 
the characteristic function of $\{x:d(x) <\eps  \}$. We have
\begin{eqnarray*}
	 N(\eps,\lam)& = & \tr [Q_{\eps}E_{\lam}Q_{\eps}]  \\
	 & = &  \tr [(Q_{\eps}E_{\lam})E_{\lam}(E_{\lam}Q_{\eps})]  \\
	 & \leq &  \norm Q_{\eps}E_{\lam}\norm ^{2}\tr [E_{\lam}] \\
	 & \leq  & c_{0}\eps^{2+1/\alp}(\lam 
+a)^{1+1/(2\alp)} N(\lam).
\end{eqnarray*}

We finally comment that lower bounds on integrals associated 
with the spectral density have recently been obtained by 
Safarov, using a coherent state method, \cite{Saf}.

\vskip 1in 
%%%%%%%%%%%%%%%%%%%%%%%%%
{\bf Acknowledgments } I would like to thank Y Safarov and M 
Solomyak for some valuable suggestions, and also acknowledge support 
under EPSRC grant number GR/L75443.
\par
\vskip 1in
%%%%%%%%%%%%%%%%%%%%%%%

\par
%%%%%%%%%%%%%%%%%%%%%%%%%%%%%%%%%%%%%%%%
\vskip 0.3in
Department of Mathematics \newline
King's College \newline
Strand \newline
London WC2R 2LS \newline
England \\
e-mail: E.Brian.Davies@kcl.ac.uk
\vfil
\end{document}